\documentclass[11pt]{article}%
\usepackage[T1]{fontenc}
\usepackage[french]{babel}
\usepackage{amsmath}
\usepackage{amsfonts}
\usepackage{amssymb}
\usepackage{graphicx}%
\def \N {\mathrm{I\!N}}
\begin{document}
\centerline {\large {\bf ON SYMMETRIC SENSITIVITY}}

\bigskip\bigskip\bigskip\centerline {Benoît CADRE and Pierre JACOB}

\bigskip\centerline {UMR CNRS 5149, Equipe de Probabilités et
Statistique}

\centerline {Université Montpellier II, CC 051, Place E. Bataillon,}

\centerline {34095 Montpellier cedex 5, FRANCE}

\bigskip\bigskip\bigskip\noindent\textbf{Abstract} We define the concept of symmetric
sensitivity with respect to initial conditions for the endomorphisms on
Lebesgue metric spaces. The idea is that the orbits of
almost every pair of nearby initial points (for the product of the invariant
measure) of a symmetrically sensitive map may diverge from a positive quantity
independent of the initial points. We study the relationships between
symmetric sensitivity and weak mixing, symmetric sensitivity and positiveness of metric entropy and 
we compute the largest sensitivity constant.



\bigskip\bigskip\bigskip\noindent\textbf{1. Introduction}

\bigskip\noindent The concept of \textit{sensitive dependence on initial
conditions} has attracted much attention in recent years and several authors
have tried to formalize it in various ways. The phrase -sensitive dependence
on initial conditions- was first used by Ruelle (1978), to indicate some
exponential rate of divergence of orbits of nearby points. More generally, it
captures the idea that a very small change in the initial condition can cause
a big change in the trajectory. Following the pioneer work by Guckenheimer
(1979), Devaney (1989) called sensitive a self-map $T\, :\, X\to X$ on the
metric space $(X,d)$ satisfying the property : there exists $\delta>0$
such that for all $x\in X$ and all $\varepsilon>0$ there is some $y\in X$
which is within a distance $\varepsilon$ of $x$ and for some $n\geq0$,
$d(T^{n}x,T^{n}y)\geq\delta$. In the last years, several authors proposed
sufficient conditions both on $T$ and $(X,d)$ to ensure the sensitivity
property (\textit{cf.} Abraham \textit{et al}, 2002, 2004, Banks
\textit{et al}, 1992, Glasner and Weiss, 1993, Guckenheimer, 1979).
  
\bigskip\noindent However, sensitive dependence on initial conditions was
first defined in Chaos Theory to measure the divergence of orbits of nearby
points, by analogy with the \textit{butterfly effect} described by the
meteorologist Ed Lorentz (for an overview in Chaos Theory, we refer the reader
to the book by Devaney, 1989). From this point of view, the above definition of
sensitivity appears to be too weak. In order to follow the sensitivity idea
drawn by the butterfly effect, one could say that $T$ is sensitive if there
exists $\delta>0$ such that for all $x,y$ in $X$, one can find $n\geq0$ with
$d(T^{n}x,T^{n}y)\geq\delta$.
However, this property appears to be too strong because it is never satisfied
by the non injective maps, such as the archetype of a chaotic map, namely the
quadratic one $Tx=4x(1-x)$ on $X=[0,1]$.

\bigskip\noindent As an attempt to weaken the previous definition of
sensitivity, we make use of tools from Ergodic Theory. From now on, we consider an endomorphism $T$ on a
probability Lebesgue space $(X,\mathcal{B},\mu)$ (\textit{cf.} Petersen, 1983,
page 16) and we fix a metric $d$ on $X$. For simplicity, we assume throughout that the support of $\mu$, denoted ${\rm supp}\, \mu$, is not reduced to a single point.

\bigskip\noindent Our sensitivity property described below is easily shown to be stronger than Guckenheimer's one. 

\bigskip\noindent\textit{\textbf{Definition} The
endomorphism $T$ is said to be symmetrically
sensitive (with respect to initial conditions) if there exists
$\delta>0$ -a sensitivity constant- such that for $\mu^{\otimes2}$-a.e. $(x,y)\in X^{2}$, one can find
$n\geq0$ with $d(T^{n}x,T^{n}y)\geq\delta$.
}

\bigskip 
\noindent Equivalently, $T$ is symmetrically
sensitive if there exists $\delta>0$ with
\[
\mu^{\otimes2}\Big(\bigcap_{n\geq 0}{\overline{T}}^{-n}\mathcal{A}_{\delta}\Big)=0,
\] 
where, here and in the following, $\overline{T}=T\times T$ is the map on $X^2$ defined by ${\overline{T}}\,:\,(x,y)\mapsto
(Tx,Ty)$ and, for any $r>0$, $\mathcal{A}_{r}$ stands for the set :
\[
\mathcal{A}_{r}:=\{(x,y)\in X^{2}\,:\,d(x,y)<r\}.
\]

\bigskip
\noindent This kind of sensitivity can be viewed as well as a property of
$\overline{T}$ which justifies in a sense the adjective {\it symmetrical}. Moreover, we note that this sensitivity property is the 
exact measure theoretic equivalent of the concept studied in Akin and Kolyada (2003) in a topological dynamic context.

\bigskip\noindent Section 2 is devoted to the computation of the sensitivity constant and 
to the case where $T$ is weakly mixing. The case where $T$ is of positive metric entropy is 
studied in Section 3.

\bigskip\noindent\textbf{2. Symmetric sensitivity, weak mixing and the sensitivity constant}

\bigskip
\noindent Observe that if $\delta$ is a sensitivity constant for $T$, then so is any positive 
$\delta'\leq \delta$. This leads to consider the following quantity, denoted $\Delta(T)$ :
$$\Delta(T)=\sup\big\{\delta\, :\, \delta \ {\rm is} \ {\rm a} \ {\rm sensitivity} \ {\rm constant} \ {\rm for} \ T\big\}.$$

\bigskip
\noindent From now on, ${\rm diam}\,(A)$ stands for the diameter of $A\subset X$ and for all $z\in X$, $r>0$, $B(z,r)$ is the open ball :
$$B(z,r)=\{x\in X\, :\, d(z,x)<r\}.$$

\bigskip
\noindent{\it {\bf Theorem 2.1} Assume that $T$ is symmetrically sensitive. Then,

\noindent (i) There exists $\delta >0$ such that for $\mu^{\otimes 2}$-a.e. $(x,y)\in X^2$, one can find a sequence $(n_k)_{k\geq 0}$ with $d(T^{n_k}x,T^{n_k}y)\geq \delta$ for all $k\geq 0$;

\noindent (ii) For $\mu^{\otimes 2}$-a.e. $(x,y)\in X^2$, one has $\sup_{n\geq 0} d(T^nx,T^ny)\geq \Delta(T)$;

\noindent (iii) $\Delta(T)\leq {\rm diam}\,  ({\rm supp}\,  \mu)$.
}

\bigskip
\noindent We first need a lemma. Notice that, since the support of $\mu$ is not reduced to a single point, there exists $\delta >0$ with $\mu^{\otimes 2}({\cal A}_{\delta})<1$. Hence, the quantity
$$a(\mu):=\sup\big\{\delta\, : \, \mu^{\otimes 2}({\cal A}_{\delta})<1\big\}$$
is positive. 

\bigskip
\noindent {\it {\bf Lemma 2.1} One has
$$a(\mu)={\rm diam}\,  ({\rm supp}\, \mu).$$
}

\bigskip
\noindent {\bf Proof} First notice that $a(\mu)\leq D:={\rm diam}\,  ({\rm supp}\, \mu)$ because for all $\varepsilon >0$, 
$\mu^{\otimes 2}({\cal A}_{D+\varepsilon})=1.$ 
Moreover, let for all $\varepsilon >0$, 
$$F_{\varepsilon}=\Big\{ (x,y)\in {\rm supp} \, \mu^{\otimes 2} :\, d(x,y)\leq D-\frac{\varepsilon}{2}\Big\},$$
where ${\rm supp} \, \mu^{\otimes 2}$ denotes the support of $\mu^{\otimes 2}$. Then, 
$\mu^{\otimes 2}(F_{\varepsilon})<1$
because  $F_{\varepsilon}$ is a closed set and 
$$F_{\varepsilon}\subsetneq  {\rm supp} \, \mu^{\otimes 2}={\rm supp} \, \mu \, \times {\rm supp} \, \mu.$$
Since ${\cal A}_{D-\varepsilon}\subset F_{\varepsilon}$, one deduces that 
$\mu^{\otimes 2}({\cal A}_{D-\varepsilon})<1$
 and hence, that $a(\mu)\geq D$ $\bullet$

\bigskip
\noindent {\bf Proof of Theorem 2.1} {\it (i)} It is a straighforward consequence of Halmos Recurence Theorem (\textit{cf.} Petersen, 1983, page 39).

\noindent {\it (ii)} If $T$ is symmetrically sensitive, then for all $\varepsilon>0$ small enough :
$$\mu^{\otimes 2} \Big( \bigcap_{n\geq 0} {\overline T}^{-n} {\cal A}_{\Delta (T)-\varepsilon} \Big)=0.$$
We get from a monotonicity argument that :
$$ \mu^{\otimes 2} \Big(\bigcup_{\varepsilon >0}  \bigcap_{n\geq 0} {\overline T}^{-n} {\cal A}_{\Delta(T)-\varepsilon} \Big)
=\lim_{\epsilon\searrow 0}\, \mu^{\otimes 2} \Big( \bigcap_{n\geq 0} {\overline T}^{-n} {\cal A}_{\Delta(T)-\varepsilon} \Big)=0,$$
hence Assertion {\it (ii)}, because
$$\bigcup_{\varepsilon >0}  \bigcap_{n\geq 0} {\overline T}^{-n} {\cal A}_{\Delta(T)-\varepsilon}=\Big\{(x,y)\in X^2\, :\, \sup_{n\geq 0} d(T^nx,T^ny)<\Delta(T)\Big\}.$$

\noindent {\it (iii)} Assume that $a(\mu) <\Delta(T)$. For any $r\in ]a(\mu),\Delta(T)[$, one has simultaneously :
$$\mu^{\otimes 2}({\cal A}_r)=1 \ {\rm and} \ \mu^{\otimes 2}\Big(\bigcap_{n\geq 0}{\overline{T}}^{-n}\mathcal{A}_r\Big)=0,$$
which is a contradiction. Therefore, $a(\mu) \geq \Delta(T)$ and {\it (iii)} is now straightforward from Lemma 2.1 $\bullet$

\bigskip\noindent\textit{\textbf{Theorem 2.2} Assume that $T$ is weakly mixing. Then, $T$ is symmetrically sensitive and moreover :
$$\Delta (T)={\rm diam} \, ({\rm supp}\, \mu).$$}

\vfill\eject
\noindent {\bf Remarks}
\begin{itemize}
\item [$\bullet$] Ergodicity is not strong enough in order to ensure the symmetric sensitivity property. Indeed, consider the case $X=I\!\!R/Z\!\!\!Z$, $\mu$ theHaar-Lebesgue measure and $d$ the natural metric on $X$. The self-map $T$
defined by $Tx=x+\theta\,(\mathrm{mod}\,1)$, where $\theta$ is an irrational
number, being an isometry for $d$, can not be symmetrically sensitive.
However, it is known to be ergodic.
\item [$\bullet$] In the case of a Guckenheimer's type definition of sensitivity, Abraham {\it et al} (2002, 2004) also provide some bounds for the largest sensitivity constant.
\item [$\bullet$]For the classical dynamical systems such as $r$-adic maps, tent map or quadratic map, one has therefore $\Delta (T)=1$. Hence, the orbits of almost all pair of nearby initial points may diverge from a quantity which is closer to 1. 
\end{itemize}

\bigskip\noindent\textbf{Proof of Theorem 2.2} Let $\delta<{\rm diam} \, ({\rm supp}\, \mu)$. Since
$\overline{T}$ is an ergodic endomorphism on $(X^{2},\mathcal{B}%
\otimes\mathcal{B},\mu^{\otimes2})$ (\textit{cf.} Petersen, 1983, page 65) and
\[
\bigcap_{n\geq0}{\overline{T}}^{-n}\mathcal{A}_{\delta}%
\]
is a $\overline{T}$-invariant set, one has
\[
\mu^{\otimes2}\Big(\bigcap_{n\geq0}{\overline{T}}^{-n}\mathcal{A}_{\delta}%
\Big)=0,
\]
because $\mu^{\otimes2}(\mathcal{A}_{\delta}^{d})<1$ by Lemma 2.1. Hence, $T$ is
symmetrically sensitive and $\Delta(T)\geq {\rm diam} \, ({\rm supp}\, \mu)$. Apply now 
Theorem 2.1 {\it (iii)}, and the theorem is proved $\bullet$

\bigskip
\noindent Observe now that for $\mu^{\otimes 2}$-a.e. $(x,y)\in X^2$, one has 
$$\sup_{n\geq 0} d(T^nx,T^ny)\leq {\rm diam} \, ({\rm supp} \, \mu).$$
Corollary 2.1 below is then a straightforward consequence of Theorem 2.1 {\it (ii)} and Theorem 2.2.

\bigskip
\noindent \textit{\textbf{Corollary 2.1} If $T$ is weakly mixing, then  for $\mu^{\otimes 2}$-a.e. $(x,y)\in X^2$ :
$$\sup_{n\geq 0} d(T^nx,T^ny)={\rm diam} \, ({\rm supp} \, \mu).$$
}

\bigskip
\noindent \textbf{3. Symmetric sensitivity and metric entropy}

\bigskip\noindent For any mesurable countable partition $\alpha$ of $X$, we
denote by $h(T, \alpha)$ the metric entropy of the transformation $T$ with
respect to the partition $\alpha$ (\textit{cf.} Petersen, 1983, Chapter 5).
Whatever being the chosen definition of sensitivity, it is usually expected
that positiveness of the entropy implies the sensitivity property ({\it cf.} Glasner and Weiss, 1993, Blanchard {\it et al}, 2002, and Abraham {\it et al}, 2004, in which positiveness of the Lyapunov exponent is considered).
Theorem 4.1 below gives an answer to this problem in the case of symmetric sensitivity. 

\bigskip\noindent\textit{\textbf{Theorem 3.1} Assume
that $T$ is ergodic and consider a finite mesurable partition $\alpha
=\{P_{1},\cdots,P_{l}\}$ of $X$. If $P_{1},\cdots,P_{l}$ are $\mu$-continuity sets for $d$ and if $h(T,\alpha)>0$, then $T$ is
symmetrically sensitive.}

\bigskip
\noindent A very similar conclusion is obtained in Blanchard {\it et al} (2002), but these authors consider the case where $T$ is a  homeomorphism on a compact space.

\bigskip\noindent\textbf{Proof} Without loss of generality, we can assume that
$h(T,\alpha)<\infty$. For all $D\in\mathcal{B}$ and $\varepsilon>0$, denote by
$D^{-\varepsilon}$ the internal $\varepsilon-$boundary of $\mathcal{D}$:%
\[
D^{-\varepsilon}=\big\{x\in D\,:\,d(x,D^{c})<\varepsilon\big\},
\]
and moreover :
\[
K_{\varepsilon}=\exp\Big(2l\sum_{i=1}^{l}\mu(P_{i}^{-\varepsilon})\Big).
\]
Since the $P_{i}$'s are $\mu$-continuity sets, $K_{\varepsilon}\rightarrow1$ as
$\varepsilon\rightarrow0$. Hence, one can choose $\delta>0$ such that
\[
K_{\delta}\,2^{-h(T,\alpha)/2}<1.\qquad(4.1)
\]
The map
\[
x\mapsto\mu\Big(\bigcap_{n\geq0}T^{-n}B(T^{n}x,\delta)\Big)
\]
defined on $X$ is $T$-invariant and moreover, according to the Fubini
Theorem,
\[
\mu^{\otimes2}\Big(\bigcap_{n\geq0}{\overline{T}}^{-n}\mathcal{A}_{\delta}\Big)=
\int_{X}\mu\Big(\bigcap_{n\geq0}T^{-n}B(T^{n}x,\delta
)\Big)\mu(dx).
\]
Consequently, by ergodicity of $T$, we have for $\mu$-a.e. $x\in X$:
\[
\mu^{\otimes2}\Big(\bigcap_{n\geq0}{\overline{T}}^{-n}\mathcal{A}_{\delta}\Big)=
\mu\Big(\bigcap_{n\geq0}T^{-n}B(T^{n}x,\delta)\Big).\qquad(4.2)
\]
We deduce from the von Neumann Ergodic Theorem that for $\mu$-a.e. $x\in X$:
\[
\frac{1}{n}\,\mathrm{card}\,\Big(k\in\{0,\cdots,n\}\,:\,T^{k}x\in\bigcup
_{i=1}^{l}P_{i}^{-\delta}\Big)\rightarrow\sum_{i=1}^{l}\mu\big(P_{i}^{-\delta
}).\qquad(4.3)
\]
We now fix a point $x\in X$ which satisfies both (4.2) and (4.3). It is
associated with it a sequence $(i_{n})_{n\geq0}\in\{1,\cdots,l\}^{\N}$ such
that $T^{n}x\in P_{i_{n}}$ for all $n\geq0$. For all $n\geq0$, we let:
\[
\mathcal{Q}_{n}=\big\{k\in\{0,\cdots,n\}\,:\,B(T^{k}x,\delta
)\ \mathrm{\nsubseteqq}\ P_{i_{k}}\big\}.
\]
For all $n\geq0$,
\begin{align*}
\mathrm{card}\,\mathcal{Q}_{n} &  \leq\mathrm{card}\,\Big(k\in\{0,\cdots
,n\}\,:\,T^{k}x\in P_{i_{k}}^{-\delta}\Big)\\
&  \leq\mathrm{card}\,\Big(k\in\{0,\cdots,n\}\,:\,T^{k}x\in\bigcup_{i=1}%
^{l}P_{i}^{-\delta}\Big)
\end{align*}
By (4.3), there exists $N_{1}\geq0$ such that for all $n\geq N_{1}$:
\[
\mathrm{card}\,\mathcal{Q}_{n}\leq2n\sum_{i=1}^{l}\mu(P_{i}^{-\delta}).
\]
For all $n\geq0$, we let:
\[
\mathcal{S}_{n}=\Big\{(s_{k})_{k=0,\cdots,n}\,:\,s_{k}\in\{1,\cdots
,l\}\ \mathrm{if}\ k\in\mathcal{Q}_{n}\ \mathrm{and}\ s_{k}=i_{k}%
\ \mathrm{if}\ k\notin\mathcal{Q}_{n}\Big\}.
\]
This set satisfies, for all $n\geq N_{1}$:
\[
\mathrm{card}\,\mathcal{S}_{n}=\exp\Big(\mathrm{card}\,\mathcal{Q}_{n}\log
l\Big)\leq\exp\Big(2nl\sum_{i=1}^{l}\mu(P_{i}^{-\delta})\Big)=K_{\delta}%
^{n}.\qquad(4.4)
\]
Now denote, for $n\geq0$ and $s\in\mathcal{S}_{n}$:
\[
L_{n,s}=\bigcap_{k=0}^{n}T^{-k}P_{s_{k}}.
\]
For all $n\geq0$, we have the following inclusions:
\begin{align*}
\bigcap_{k=0}^{n}T^{-k}B(T^{k}x,\delta) &  =\bigcap_{k\in\mathcal{Q}%
_{n}^{c}\cup\mathcal{Q}_{n}}T^{-k}B(T^{k}x,\delta)\\
&  \subset\bigcap_{k\in\mathcal{Q}_{n}^{c}}T^{-k}P_{i_{k}}\\
&  \subset\bigcup_{s\in\mathcal{S}_{n}}L_{n,s}.\qquad(4.5)
\end{align*}
Now fix $\varepsilon\in]0,h(T,\alpha)/2[$. By the Entropy Equipartition
Property (\textit{cf.} Petersen, 1983, page 263), there exists $N_{2}\geq0$ such
that for all $n\geq N_{2}$, the elements of $\{1,\cdots,l\}^{n+1}$ can be
divided into two disjoints classes, $\mathcal{G}_{n}$ and $\mathcal{B}_{n}$,
such that :
\[
\mu\Big(\bigcup_{s\in\mathcal{B}_{n}}L_{n,s}\Big)\leq\varepsilon,
\]
and, for all $s\in\mathcal{G}_{n}$ :
\[
\mu(L_{n,s})\leq2^{-n(h(T,\alpha)-\varepsilon)}.
\]
We then deduce from (4.4) and (4.5) that for all $n\geq\max(N_{1},N_{2})$:
\begin{align*}
\mu\Big(\bigcap_{k=0}^{n}T^{-k}B(T^{k}x,\delta)\Big) &  \leq
\mu\Big(\bigcup_{s\in\mathcal{S}_{n}\cap\mathcal{B}_{n}}L_{n,s}\Big)+\sum
_{s\in\mathcal{S}_{n}\cap\mathcal{G}_{n}}\mu(L_{n,s})\\
&  \leq\varepsilon+\mathrm{card}\,\mathcal{S}_{n}\,\max_{s\in\mathcal{G}_{n}%
}\,\mu(L_{n,s})\\
&  \leq\varepsilon+K_{\delta}^{n}\,2^{-n(h(T,\alpha)-\varepsilon)}\\
&  \leq\varepsilon+K_{\delta}^{n}\,2^{-nh(T,\alpha)/2},
\end{align*}
where the latter inequality comes from the fact that $\varepsilon
<h(T,\alpha)/2$. Letting $n\rightarrow\infty$, we deduce from (4.1) that for
all $\varepsilon>0$ small enough :
\[
\lim_{n}\mu\Big(\bigcap_{k=0}^{n}T^{-k}B(T^{k}x,\delta)\Big)\leq
\varepsilon,
\]
hence, letting $\varepsilon\rightarrow0$ :
\[
\mu\Big(\bigcap_{n\geq0}T^{-n}B(T^{n}x,\delta)\Big)=0.
\]
Finally, we deduce from (4.2) and the choice of $x$ that $T$ is symmetrically
sensitive $\bullet$

\bigskip\bigskip\centerline {\bf REFERENCES}

\bigskip\bigskip

\noindent Abraham, C., Biau, G. and Cadre, B. (2002). Chaotic Properties of
Mappings Defined on a Probability Space, \textit{J. Math. Anal. Appl.}
\textbf{266}, pp. 420-431.

\bigskip\noindent  Abraham, C., Biau, G. and Cadre, B. (2004). On Lyapunov
Exponent and Sensitivity, \textit{J. Math. Anal. Appl.} \textbf{290}, pp. 395-404.

\bigskip\noindent Akin, E. and Kolyada, S. (2003). Li-Yorke sensitivity, \textit{Nonlinearity} \textbf{16}, pp. 1421-1433.

\bigskip\noindent Banks, J., Brooks, J., Cairns, G., Davis, G. and Stacey,
P. (1992). On Devaney's Definition of Chaos, \textit{Am. Math. Mon.} \textbf{99}, pp. 332-334.

\bigskip\noindent Blanchard, F., Host, B., Ruette, S. (2002). Asymptotic Pairs in Positive-Entropy Systems, \textit{Ergodic Theory Dynam. Systems} \textbf{22}, pp. 671-686.

\bigskip\noindent Devaney, R.L. (1989). \textit{Chaotic Dynamical Systems}, 2nd edn (New-York: Addison-Wesley).

\bigskip\noindent Glasner, E. and Weiss, B. (1993). Sensitive Dependence on
Initial Conditions, \textit{Nonlinearity} \textbf{6}, pp. 1067-1075.

\bigskip\noindent  Guckenheimer, J. (1979). Sensitive Dependence to Initial
Conditions for One Dimensional Maps, \textit{Commun. Math. Phys.} \textbf{70},
pp. 133-160.

\bigskip\noindent Petersen, K. (1983). \textit{Ergodic Theory}, (Cambridge: Cambridge University Press).

\bigskip\noindent Ruelle, D. (1978). Dynamical Systems with
Turbulent Behavior, {\it Mathematical Problems in Theoretical Physics} (Lectures Notes in Physics), (Berlin: Springer).

\end{document}